\documentclass[12pt]{article}
\usepackage{amsmath,amssymb}
\input{epsf}

\newtheorem{theorem}{Theorem}

\newtheorem{proposition}{Proposition}

\newtheorem{question}{Question}

\newenvironment{proof}{{\bf Proof.}}{\hfill\rule{2mm}{2mm}}

\newtheorem{remarka}{Remark}
\newenvironment{remark}{\begin{remarka}\rm}{\end{remarka}}

\def\sgn {{\rm sgn}}

\title{\bf Perfect dominating sets in the Cartesian products of prime cycles.}
\author{{\bf Hamed Hatami$^a$ and Pooya Hatami$^b$} \\
$^a${\small\it Department of Computer Science}\\
{\small {University of Toronto}}\\
$^b${\small\it Department of Mathematical Sciences}\\
{\small Sharif University of Technology}}

\date{}

\begin{document}
\maketitle

\begin{abstract}
We study the structure of a minimum dominating set of
$C_{2n+1}^n$, the Cartesian product of $n$ copies of the cycle of
size $2n+1$, where $2n+1$ is a prime.
\end{abstract}


\section{Introduction}        
Let $G$ and $H$ be two graphs. The \emph{Cartesian product} of $G$
and $H$ is a graph with vertices $\{(x,y): x \in G, y \in H\}$
where $(x,y) \sim (x',y')$ if and only if $x=x'$ and $y \sim y'$,
or $x \sim x'$ and $y=y'$. Let $G^n$ denote the Cartesian product
of $n$ copies of $G$. This article deals with $C_{2n+1}^n$ where
$C_{2n+1}$ is the cycle of size $p:=2n+1$ and $p$ is a prime.

For our purpose, it is more convenient to view the vertices of
$C_{2n+1}^n$ as the elements of the group
$G:=\mathbb{Z}_{2n+1}^n$. Then $x \sim y$ if and only if $x-y=\pm
e_i$ for some $i\in [n]$, where $e_i=(0,\ldots,1,\ldots,0)$ is the
unit vector with $1$ at the $i$th coordinate. In other words,
$C_{2n+1}^n$ is the Cayley graph
$\Gamma(\mathbb{Z}_{2n+1}^n,{\cal U})$ over the group
$\mathbb{Z}_{2n+1}^n$ with the set of generators ${\cal U}=\{\pm
e_1,\ldots,\pm e_n\}$. From this point on, to emphasis the group
structure of the graph we will use the Cayley graph notation
$\Gamma(\mathbb{Z}_{2n+1}^n,{\cal U})$ instead of the Cartesian
product notation of $C_{2n+1}^n$.

Let $u$ and $v$ be two vertices of a graph $G$. We say that $u$
dominates $v$ if $u=v$ or $u \sim v$. A subset $S$ of the vertices
of $G$ is called a \emph{dominating set}, if every vertex of $G$
is dominated by at least one vertex of $S$. A dominating set is
\emph{perfect}, if no vertex is dominated by more than one vertex.

\begin{remark}
\label{rmk:lalgebra1}  Let $G$ be a graph with $m$ vertices.
Every function $f:V(G) \rightarrow \mathbb{C}$ can be viewed as a
vector $\vec{f} \in \mathbb{C}^m$. Let $A$ denote the adjacency
matrix of $G$. Note that $f:V(G)\rightarrow \{0,1\}$ is the
characteristic function of a perfect dominating set if and only if
$(A+I)\vec{f}=\vec{1}$.
\end{remark}

We are interested in perfect dominating sets of
$\Gamma(\mathbb{Z}_{2n+1}^n,{\cal U})$. Note that for an
$r$-regular graph $G=(V,E)$ a dominating set is perfect if and
only if it is of size $|V|/(r+1)$. Since
$\Gamma(\mathbb{Z}_{2n+1}^n,{\cal U})$ is $2n$-regular and has
$(2n+1)^n$ vertices, a dominating set is perfect if and only if
it is of size $(2n+1)^{n-1}$.

 Fix an arbitrary
$(\epsilon_1,\ldots,\epsilon_{n-1}) \in \{-1,1\}^{n-1}$, and a $k
\in \{0,\ldots,2n\}$. The set
\begin{equation}
\label{eq:linear} \{(x_1,\ldots,x_{n-1},k+\sum_{i=1}^{n-1}
\epsilon_{i}(i+1)x_i): x_i \in \mathbb{Z}_{2n+1} \ \forall i \in
[n-1]\}
\end{equation}
forms a perfect dominating set in $\Gamma(\mathbb{Z}_{2n+1}^n,{\cal
U})$, where the additions are in $\mathbb{Z}_{2n+1}$. To see this
consider $y=(y_1,\ldots,y_n) \in \mathbb{Z}_{2n+1}^n$. Let
$t=k+\sum_{i=1}^{n-1} \epsilon_{i}(i+1)y_i$, and  $\Delta = t-y_n
\mod 2n+1$ so that $|\Delta| \le n$. If $\Delta \in \{-1,0,1\}$ then
$y$ is dominated by $(y_1,\ldots,y_{n-1},t)$. If $\Delta \not\in
\{-1,0,1\}$, then with the notation $j:=|\Delta|-1$, $y$ is adjacent
to $(y_1,\ldots,y_{j-1},y_j-\epsilon_j \times
\sgn(\Delta),y_{j+1},\ldots,y_n)$, which can  easily be seen that is
in the considered set.

There are many results in the direction of constructing  perfect
dominating sets in the Cartesian product of cycles
(see~\cite{MR1140510} and its references). However the authors are
unaware of any result in the direction of characterizing the
structure of perfect dominating sets. We consider the simplest case
$\Gamma(\mathbb{Z}_{2n+1}^n,{\cal U})$ where $2n+1$ is a prime. Even
in this simple case we are unable to characterize all the perfect
dominating sets. However we prove the following theorem in this
direction.

\begin{theorem}
\label{thm:main} Let $2n+1$ be a prime and $S \subseteq
\Gamma(\mathbb{Z}_{2n+1}^n,{\cal U})$ be a perfect dominating
set. Then for every $(x_1,\ldots,x_n) \in \mathbb{Z}_{2n+1}^n$
and every $i \in [n]$,
$$ |S \cap \{(y_1,\ldots,y_n) : y_j=x_j \ \forall j \neq i\}|=1.$$
\end{theorem}

Theorem~\ref{thm:main} says that when $2n+1$ is a prime, every
parallel-axis line contains exactly one point from every perfect
dominating set of $\Gamma(\mathbb{Z}_{2n+1}^n,{\cal U})$. It is
easy to construct examples to show that the condition of $2n+1$
being a prime is necessary~\cite{Golomb}.

Let ${\cal F}$ be a family of sets. For $S \in {\cal F}$, a set $D
\subseteq S$ is called a \emph{defining set} for $(S, {\cal F})$
(or for $S$ when there is no ambiguity), if and only if $S$ is the
only superset of $D$ in ${\cal F}$. The size of the minimum
defining set for $(S, {\cal F})$ is called its defining number.
Defining sets are studied for various families of ${\cal F}$
(See~\cite{MR2011736} for a survey on the topic). Let ${\cal F}$
be the family of all minimum dominating sets of
$\Gamma(\mathbb{Z}_{2n+1}^n,{\cal U})$. Note that since
$\Gamma(\mathbb{Z}_{2n+1}^n,{\cal U})$ is regular and contains at
least one perfect dominating set, a set $S \subseteq V(G)$ is a
minimum dominating set if and only if it is a perfect dominating
set. In~\cite{MR1480799} Chartrand et al. studied the size of
defining sets of ${\cal F}$ for $n=2$. Based on this case they
conjectured that the smallest defining set over all minimum
dominating sets of $\Gamma(\mathbb{Z}_{2n+1}^n,{\cal U})$ is of
size exactly $n$. As it is noticed by Richard Bean~[private
communication], the conjecture fails for $n=3$, as in this case
there are perfect dominating sets with defining number $2$ (See
Remark~\ref{rem:case7}). So far there is no nontrivial bound
known for the defining numbers of minimum dominating sets of
$\Gamma(\mathbb{Z}_{2n+1}^n,{\cal U})$. We prove the following
theorem.

\begin{theorem}\label{thm:defining1}
Let $2n+1$ be a prime and ${\cal F}$ be the family of all minimum
dominating sets of $\Gamma(\mathbb{Z}_{2n+1}^n,{\cal U})$. Every
$S \in {\cal F}$ has a defining set of size at most $n! 2^n$.
\end{theorem}

The proof of Theorem~\ref{thm:main} uses Fourier analysis on
finite Abelian groups. In Section~2 we review Fourier analysis on
$\mathbb{Z}_p^n$. Section~3 is devoted to the proof of
Theorem~\ref{thm:defining1}. Section~4 contains further
discussions about the defining sets of minimum dominating sets of
$\Gamma(\mathbb{Z}_{2n+1}^n,{\cal U})$.

\section{Background}
In this section we introduce some notations and review Fourier
analysis on $G=\mathbb{Z}_p^n$. For a nice and more detailed, but
yet brief introduction we refer the reader to~\cite{Alon}. See
also~\cite{Groups} for a more comprehensive reference.

 Aside from its group structure we will also
think of $G$ as a measure space with the uniform (product)
measure, which we denote by $\mu$. For any function $f:G
\rightarrow \mathbb{C}$, let
$$\int_G f(x)dx=\frac{1}{|G|}\sum_{x \in G} f(x).$$
The inner product between two functions $f$ and $g$ is
$\langle f,g \rangle=\int_G f(x) \overline{g(x)} dx.$
Let
$$\omega=e^{2\pi i/p},$$
where $i$ is the imaginary number. For any $x \in G$, let
$\chi_x:G \rightarrow \mathbb{C}$ be defined as
$$\chi_x(y)=\omega^{\sum_{i=1}^n x_i y_i}.$$

It is easy to see that these functions form an orthonormal basis.
So every function $f:G \rightarrow \mathbb{C}$ has a unique
expansion of the form $f=\sum \widehat{f}(x)\chi_x,$ where
$\widehat{f}(x)=\langle f,\chi_x\rangle$ is a complex number.

\section{Proof of Theorem~\ref{thm:main}}
Let $\vec{0}=(0,\ldots,0)$, $\vec{1}=(1,\ldots,1)$, and
$e_i=(0,\ldots,1,\ldots,0)$, the unit vector with $1$ at the
$i$-th coordinate. Let $p=2n+1$ be a prime and $S$ be a perfect
dominating set in $G$, and let $f$ be the characteristic function
of $S$, i.e. $f(x)=1$ if $x\in S$ and $f(x)=0$ otherwise. Let
$$D=\{\pm e_1, \pm e_2, \ldots, \pm e_n\},$$
be the set of unit vectors and their negations. For every $\tau
\in D$ define $f_\tau(x)=f(x+\tau)$. Note that
$$ \widehat{f}_\tau(y)= \int f(x + \tau) \overline{\chi_y(x)} dx= \int
f(x) \overline{\chi_y(x-\tau)}dx=\int f(x)
\overline{\chi_y(x)}\chi_y(\tau) dx = \widehat{f}(y)
\chi_y(\tau).$$
 Let
$$g= f + \sum_{\tau \in D} f_\tau.$$
We have

\begin{equation}
\label{eq:exp-g} g= \left(\sum_{y \in G}
\widehat{f}(y)\chi_y\right) + \sum_{\tau \in D} \sum_{y \in G}
\widehat{f_\tau}(y) \chi_y = \sum_{y \in G} \widehat{f}(y)
\left(\sum_{\tau \in D \cup \{\vec{0}\}} \chi_y(\tau)\right)
\chi_y.
\end{equation}
Since $f$ is the characteristic function of a perfect dominating
set, we have $g(x)=1$, for every $x \in G$. So $g =
\chi_{\vec{0}}$. By uniqueness of Fourier expansion, for every $y
\neq \vec{0}$,
\begin{equation}
\label{eq:expanded} 0=\widehat{g}(y)=\widehat{f}(y) \sum_{\tau \in
D \cup \{\vec{0}\}} \chi_y(\tau)=
\widehat{f}(y)\left(1+\sum_{i=1}^n \omega^{y_i}+\sum_{i=1}^n
\omega^{-y_i}\right).
\end{equation}
Now we turn to the key step of the proof. Since $2n+1$ is a prime,
(\ref{eq:expanded}) implies that whenever $\widehat{f}(y) \neq 0$,
we have
\begin{equation}
\label{eq:nonz} \{y_1,\ldots,y_n\} \cup
\{-y_1,\ldots,-y_n\}=\{1,\ldots,2n\}.
\end{equation}
Denote the set of all $y$ satisfying~(\ref{eq:nonz}) by ${\cal
Y}$. For $1 \le i \le n$, let
$$D_i=\{k e_i: 0 \le k \le 2n\}.$$
Define $g_i=\sum_{\tau \in D_i} f_\tau$. Similar to
(\ref{eq:exp-g}), we get
$$g_i=\sum  \widehat{f}(y) \left(\sum_{\tau
\in D_i} \chi_y(\tau)\right) \chi_y.$$
When $y \in {\cal Y}$, since $y_i \neq 0$, we have
$$\sum_{\tau \in D_i} \chi_y(\tau)=  \sum_{k=0}^{2n}
\omega^{k y_i}=0.$$
When $y {\not\in} {\cal Y}$ and $y \neq \vec{0}$,
$\widehat{f}(y)=0$. So

\begin{equation}
\label{lines} g_i=\left(\widehat{f}(0) \sum_{\tau \in D_i}
\chi_{\vec{0}}(\tau)\right) \chi_{\vec{0}}=\chi_{\vec{0}}=1.
\end{equation}
Note that $g_i(x)$ counts the number of elements in $S \cap
\{(y_1,\ldots,y_n) : y_j=x_j \ \forall j \neq i\}$. This completes
the proof.

\begin{remark}
\label{rmk:lalgebra} The above proof can be translated to the
language of linear algebra (However in the linear algebra language
the key observation~(\ref{eq:nonz}) becomes less obvious). Indeed,
let $m=(2n+1)^n$ denote the number of vertices. From
Remark~\ref{rmk:lalgebra1} we know that
$f:\mathbb{Z}_{2n+1}^n\rightarrow \mathbb{C}$ is the
characteristic function of a perfect dominating set if and only if
$(A+I)\vec{f}=\vec{1}$, where $A$ is the adjacency matrix of
$\Gamma(\mathbb{Z}_{2n+1}^n,{\cal U})$. The reader may notice that
in the proof of Theorem~\ref{thm:main}, $\vec{g}=(A+I)\vec{f}$,
and thus (\ref{eq:exp-g}) shows that $\vec{\chi_y}$ is a family of
orthonormal eigenvectors of $A+I$. Moreover, among these
eigenvectors, the ones that correspond to the $0$ eigenvalue are
exactly $\vec{\chi_y}$ with $y \in {\cal Y}$. Hence the rank of
$A+I$ is $m-|{\cal Y}|=(2n+1)^n-n! 2^n$. We will use this fact in
the proof of Theorem~\ref{thm:defining1}.
\end{remark}

\section{Proof of Theorem~\ref{thm:defining1}}
As it is observed in Remark~\ref{rmk:lalgebra1}, every perfect
dominating set of $\Gamma(\mathbb{Z}_{2n+1}^n,{\cal U})$
corresponds to a zero-one vector $\vec{f} \in \mathbb{C}^m$ that
satisfies $(A+I)\vec{f}=\vec{1}$. Let
$$V:={\rm span}\{\vec{f}: f \in {\cal F}\}.$$
Trivially
$$\dim V \le 1+(m-{\rm rank}(A+I)) = 1+n! 2^n.$$
Also for a subset $D$ of vertices of
$\Gamma(\mathbb{Z}_{2n+1}^n,{\cal U})$, define
$$V_D:={\rm span}\{\vec{f}: \mbox{$f \in {\cal F}$ and $\forall x \in D,\ f(x)=1$}\},$$
 Note that
$V=V_\emptyset$.

To prove Theorem~\ref{thm:defining1} we start from $D=\emptyset$.
At every step, if $D$ does not extend uniquely to $S$, then there
exists a vertex $v \in S$ such that $\dim V_{D \cup \{v\}} < \dim
V_{D}$; we add $v$ to $D$. Since $\dim V_\emptyset \le 1+n! 2^n$,
we can obtain a set $D$ of size at most $n! 2^n$ such that the
dimension of $V_D$ is at most $1$. This completes the proof as
there is at most one non-zero, zero-one vector in a vector space
of dimension $1$.

\section{Future directions}

We ask the following question:

\begin{question}
\label{ques:linear} For a prime $2n+1$, are there examples of
perfect dominating sets in $\Gamma(\mathbb{Z}_{2n+1}^n,{\cal U})$
that are not of the form (\ref{eq:linear})?
\end{question}

If the answer to Question~\ref{ques:linear} turns out to be
negative, then we can improve the bound of
Theorem~\ref{thm:defining1}:

\begin{proposition}\label{prop} Let $p=2n+1$ be a prime, and let ${\cal T}$ denote the set of perfect
dominating sets of the form~(\ref{eq:linear}). Every $(S,{\cal
T})$ where $S \in {\cal T}$ has a defining set of size $1+\lceil
\frac{n-1}{\lfloor \log_{2}{p} \rfloor}\rceil$.
\end{proposition}
\begin{proof}
Suppose that $S \in {\cal T}$. Then $S$ is of the form:
$$
\label{eq:linear2} \{(x_1,\ldots,x_{n-1},k+\sum_{i=1}^{n-1}
\epsilon_{i}(i+1)x_i): x_i \in \mathbb{Z}_p \ \forall i \in
[n-1]\}.
$$
Let $m=\lfloor \log_{2}{p} \rfloor$. We will use the easy fact
that for any $c \in \mathbb{Z}_p$, the equation $\sum_{i=0}^{m-1}
\epsilon_i 2^i =_p c$ has at most one solution $(\epsilon_0,
\epsilon_1,...,\epsilon_{m-1})\in \{-1, +1\}^m$. For $i,j \ge 0$,
define $\alpha_{i,j} \in \mathbb{Z}_p$ to be the solution to
$(i+j+1) \alpha_{i,j}=_p 2^{j}$.

Let $u=(0,0,...,0,b)$ be the unique vertex in $S$ with the first
$n-1$ coordinates equal to $0$, and for every $1 \le i \le n-1$
consider the unique vector
$$u_i=(0,...,0,\alpha_{i,0},\alpha_{i,1},...,\alpha_{i,k_i},0,...,0,b_i) \in S,$$
where $\alpha_{i,0}$ is in the $i$th coordinate and
$k_i=\min(m-1,n-i-1)$. We claim that the set $D=\{ u, u_0, u_{m},
..., u_{m (\lceil\frac{n-1}{m}\rceil-1)}\}$ is a defining set for
$(S,{\cal T})$. Since $S$ is of form~(\ref{eq:linear}), clearly
$k=b$, and for every $0 \le i \le \lceil \frac{n-1}{m}\rceil-1$,
we have:

\begin{equation*}
b_{mi}-b =\sum_{j=0}^{k_{mi}}
\epsilon_{mi+j}(mi+j+1)\alpha_{mi,j}=\sum_{j=0}^{k_{mi}}
\epsilon_{mi+j} 2^j.
\end{equation*}

The above equation has only one solution
$(\epsilon_{mi},\epsilon_{mi+1},...,\epsilon_{mi+k_{mi}})\in
\{-1,+1\}^{k_{mi}+1}$. Considering this for all $u_{mi} \in D$
determines $(\epsilon_1, \epsilon_2, ..., \epsilon_{n-1})$. Thus
the set $D$ is a defining set for $(S,{\cal T})$.
\end{proof}

\begin{remark}
\label{rem:case7} For $n=2,3$ the answer to
Question~\ref{ques:linear} is negative. Thus when $n=3$,
Proposition~\ref{prop} implies that there is a defining set of
size $2$ for a perfect dominating set. This disproves the
conjecture of~\cite{MR1480799} which is already observed by
Richard Bean~[private communication].
\end{remark}

\bibliographystyle{plain}
\bibliography{dom}

\end{document}